\documentclass[conference]{IEEEtran}
\IEEEoverridecommandlockouts
\usepackage{cite}
\usepackage{amsmath}
\usepackage{amssymb}
\usepackage{amsfonts}
\usepackage{amsthm}
\usepackage[utf8]{inputenc}
\usepackage[english]{babel}
\usepackage{tikz}
\usetikzlibrary{shapes,arrows,chains}
\tikzstyle{block} = [rectangle, draw,text centered, rounded corners]
\tikzstyle{line} = [draw, -latex']
\usepackage{blindtext}

\usepackage{subcaption}
\usepackage{ulem}
\usepackage{float}

\newtheorem{theorem}{Theorem}[section]

\newtheorem{property}[theorem]{Property}

\theoremstyle{remark}

\theoremstyle{definition}

\usepackage{algorithmic}
\usepackage{graphicx}
\usepackage{textcomp}
\usepackage{xcolor}
\def\BibTeX{{\rm B\kern-.05em{\sc i\kern-.025em b}\kern-.08em
    T\kern-.1667em\lower.7ex\hbox{E}\kern-.125emX}}

\begin{document}

\title{Maximal Positive Invariant Set Determination for Transient Stability Assessment in Power Systems}

\author{\IEEEauthorblockN{Antoine Oustry}
\IEEEauthorblockA{
\textit{Ecole Polytechnique}\\
Palaiseau, France \\
antoine.oustry@polytechnique.edu}
\and
\IEEEauthorblockN{Carmen Cardozo}
\IEEEauthorblockN{Patrick Panciatici, \textit{IEEE fellow}}
\IEEEauthorblockA{\textit{RTE R\&D} \\
Versailles, France \\
carmen.cardozo@rte-france.com}
\and
\IEEEauthorblockN{Didier Henrion}
\IEEEauthorblockA{\textit{CNRS-LAAS} \\
Toulouse, France \\
\textit{and FEL-\v CVUT} \\
Prague, Czechia\\
henrion@laas.fr}}

\maketitle

\begin{abstract}
This paper assesses the transient stability of a synchronous machine connected to an infinite bus through the notion of invariant sets. The problem of computing a conservative approximation of the maximal positive invariant set is formulated as a semi-definitive program based on occupation measures and Lasserre's relaxation. An extension of the proposed method into a robust formulation allows us to handle Taylor approximation errors for non-polynomial systems. Results show the potential of this approach to limit the use of extensive time domain simulations provided that scalability issues are addressed. 
\end{abstract} 


\section{Introduction}

Although a classic definition of dynamic system stability does apply to power systems, this notion has been traditionally classified into different categories depending on the variables (rotor angle, voltage magnitude or frequency), the time scale (short and long term) of interest~\cite{{kundur1994power}}, as well as the size of the disturbance. 
In particular, transient stability refers to the ability of the power system to maintain synchronism when subjected to a severe disturbance and it focuses on the evolution of generator rotor angles over the first 
seconds that follow. 

Indeed, a short circuit at a synchronous generator's terminal reduces its output voltage, and with it, the power injected into the network. 
The received mechanical power is then stored in the rotor mass as kinetic energy producing a speed increase. If the voltage is not restored within a certain time for the specific fault, known as the Critical Clearing Time (CCT), the unit loses synchronism, i.e. the rotor angle diverges.

Transmission System Operators (TSO) are responsible for the power system security and must prevent this to happen as a consequence of any plausible N-1 situation. Hence, TSOs constantly perform intensive time-domain nonlinear simulations and take actions if needed to ensure transient stability. Historically, the simulated scenarios could be limited to a manageable set of given initial conditions and predefined faults. However, with the changing operational environment of electrical power systems these critical conditions become harder to identify. Renewable energy sources and new architectures of intraday and balancing markets add uncertainty and variability to the production plan, enlarging the set of possible initial conditions that TSOs have to consider.

Therefore, the research of new methods for assessing the transient stability of classical power systems has drawn academia and industry attention. The computation of Regions of Attraction (ROA) for this purposes appeared as an interesting idea. Indeed, a ROA provides the set of acceptable (post-fault) conditions of a dynamic system which are known to reach a given target set in a specified time. They can be obtained through the construction of polynomial Lyapunov functions~\cite{anghel2013algorithmic,kalemba2018Lypunov}, as well as using the notion of occupation measures and Lasserre's hierarchy~\cite{henrion2014convex,korda2013inner}.
As long as the dynamics of the system is polynomial, both formulations yield a moment-sum-of-squares (SOS) optimization program that can be efficiently solved by semi-definite programming (SDP), a particular class of efficient convex optimization tools.

Within the framework of occupation measures, we propose to assess the transient stability of a power system by computing Maximal Positively Invariant (MPI) sets which will simply exclude all diverging trajectories without fixing any arbitrary target set and reaching time. We consider a Single Machine Infinite Bus (SMIB) test system, based on different non-polynomial models and classical hypotheses of electromechanical analysis. The originality of this work lies on the reformulation of the problem presented in~\cite{MCIOUTER} into the inner approximations of the MPI set and its application to the transient stability study of a synchronous machine (SM). 

The main contributions of this work are:

\begin{enumerate}
    \item Formulation of the MPI set inner approximation problem for a polynomial dynamic system constrained on an algebraic set. 
    \item Its extension into a robust form that ensures the conservativeness of the MPI set for a non-polynomial system as long as its approximation error can be bounded.
    \item Computation of CCT bounds without simulation of the post-fault system, but from the evaluation of the polynomial describing the inner/outer approximations of the MPI set along the trajectory of the faulted system.
\end{enumerate}

Section~\ref{sec:Modelling} presents the polynomial reformulation of  three different SMIB models. Section~\ref{sec:Inner} describes the MPI set inner approximation method while Section~\ref{sec:Robust} includes its robust form. 
Numerical results are analyzed in Section~\ref{sec:Results}. Conclusion and future work are discussed in Section~\ref{sec:Conclusions}.

\section{Polynomial reformulation of the SM Models}
\label{sec:Modelling}

In this work we consider three different SM models. The second order model (2\textsuperscript{nd} OM) is described as follows:
\begin{equation}\begin{small}
  \left\{
  \begin{array}{l}
    \dot{\delta}(t) = \omega_n (\omega_s(t) -\omega_i) \\
    2H \dot{\omega_s}(t) = C_m -\frac{1}{\omega_s(t)} \frac{V_s V_i}{X_l} \sin(\delta(t)) - D  (\omega_s(t) -\omega_i)   \\
  \end{array}
  \right.
\label{eq:2ndOM}
\end{small}\end{equation}
where $\delta$ (radians) is the angle difference between the generator and the infinite bus, $H$ (MWs/MVA) the inertia constant, $\omega_n$ (radians/s) the nominal frequency, $\omega_s$ the generator speed, $\omega_i$ the infinite bus speed, $V_s$ the generator voltage, $V_i$ the infinite bus voltage, $C_m$ the mechanical torque, $X_l$ the line reactance and $D$ the damping factor, all in per unit (p.u.). 

The third order model (3\textsuperscript{rd} OM) takes into account the dynamics of the transient electromotive force ($e'_q$) considering a constant exciter output voltage ($E_{vf}$):
\begin{equation}
\begin{small}
  \left\{
  \begin{array}{l}
    \dot{\delta}(t) = \omega_n (\omega_s(t) -\omega_i) \\
    2H \dot{\omega_s}(t) = C_m -\frac{1}{\omega_s(t)} \frac{V_s V_i}{X_l} \sin(\delta(t)) - D  (\omega_s(t) -\omega_i)   \\
    T^{'}_{d0} \dot{e'_q}(t) = E_{vf} - e'_q(t) + \frac{x_d-x^{'}_d}{x^{'}_d+X_l}(V_i \cos{\delta(t)} - e'_q(t)) \\
  \end{array}
  \right.
\label{eq:3thOM}
\end{small}
\end{equation}
where $x_d$ and $x'_d$ are the SM steady state and transient direct axis reactances, and $T^{'}_{d0}$ is the direct axis open-circuit transient time-constant.
The fourth order model (4\textsuperscript{th} OM) includes a voltage controller:
\begin{equation}
\begin{small}
  \left\{
  \begin{array}{l}
    \dot{\delta}(t) = \omega_n (\omega_s(t) -\omega_i) \\
    2H \dot{\omega_s}(t) = C_m -\frac{1}{\omega_s(t)} \frac{V_s V_i}{X_l} \sin(\delta(t)) - D  (\omega_s(t) -\omega_i)   \\
    T^{'}_{d0} \dot{e'_q}(t) = E_{vf}(t) - e'_q(t) + \frac{x_d-x^{'}_d}{x^{'}_d+X_l}(V_i \cos{\delta(t)} - e'_q(t)) \\
      T_E \dot{E_{vf}}(t) = \kappa (V_{ref} - V_s(t)) - E_{vf}(t)\\
  \end{array}
  \right.
\label{eq:4thOM}
\end{small}
\end{equation}
where now the exciter output voltage $E_{vf}(t)$ is time varying,
\[
V_s(t)=\sqrt{(\frac{x_q V_i \sin{\delta(t)}}{x_q+X_l})^2+(\frac{x_d V_i \cos{\delta(t)} + X_l e'_q(t)}{x^{'}_d+X_l})^2},
\]
$x_q$ is the quadrature axis reactance, $V_{ref}$ is the SM reference voltage and $\kappa$ is the controller gain, all in p.u.  
These models include non-polynomial terms on $\delta$ (trigonometric function), $\omega$ (inverse function) and also $e'_q$ (square root). In the sequel we explain how to derive polynomial models by reformulations.

\subsection{Variable change for exact equivalent}

As demonstrated in~\cite{matteo2018Lypunov}, the trajectories and stability properties of the system are preserved when using the following endogenous transformation:
\begin{equation}
  \Phi := \left\{
  \begin{array}{l}
   ]-\pi,\pi[ \times  ]-\omega_{M},\omega_{M}[ \to \mathcal{C} \times ]-\omega_{M},\omega_{M}[ \\
   (\delta,\omega) \mapsto (\cos(\delta),\sin(\delta),\omega)
  \end{array}
  \right.
\label{eq:change_variables}
\end{equation}
where $\omega_M$ is an upper bound on $\omega$ and $\mathcal{C}=\{(x,y)\in \mathbb{R}^2, x^2+y^2=1, x>-1\}$. Then, the SM 2\textsuperscript{nd} OM becomes polynomial at the price of increasing the dimension of the state space and adding an algebraic constraint. 

\subsection{Taylor Approximation}

The polynomial reformulation of the inverse function and the square root, whose arguments have limited variations in the post-fault system, is achieved using a classic Taylor series expansion. Without loss of generality, we set $\omega(t)=\omega_s(t) - \omega_i$ as the speed deviation of the SM and $\omega_i=$ 1 p.u., such that:
\begin{equation}
\frac{1}{1+\omega} = 1 - \omega + \omega ^2  + o(\omega^2) \\
\label{eq:DLomega}
\end{equation}

\vspace{-3mm}

\begin{equation}
V_s = V_s^{eq} (1+\frac{h}{2}-\frac{h^2}{8})+ o(h^2)\\
 \label{eq:DLVs}
\end{equation}
where $V_s^{eq}$ is the terminal voltage at an equilibrium point and $h = [(\frac{x_q V_i y}{x_q+X_l})^2+(\frac{x_d V_i x + X_l e'_q}{x^{'}_d+X_l})^2] /(V_s^{eq})^2 -1 $. 

\subsection{Polynomial Model}
The 4\textsuperscript{th} SM model is now expressed as a polynomial model:
\begin{equation}\begin{small}
  \left\{
  \begin{array}{l}
    \dot{x}(t) = - \omega_n \omega(t)  y(t) \\
    \dot{y}(t) = \omega_n \omega(t)  x(t) \\
     2H \dot{\omega}(t) = C_m -(1 - \omega(t) + \omega ^2 (t)) \frac{V_s V_i}{X_l} y(t) - D \omega (t) \\
    T^{'}_{d0} \dot{e'_q}(t) = \Bar{E}_{vf}(t) - e'_q(t) + \frac{x_d-x^{'}_d}{x^{'}_d+X_l}(V_i \cos{\delta}(t) - e'_q(t)) \\
      T_E \dot{E_{vf}}(t) = \kappa (V_{ref} - (V_s^{eq} (1+0.5h(t)-0.125h(t)^2)) - E_{vf}(t)\\
    x(t)^2+y(t)^2 = 1, x(t) > -1, \omega(t) \in  ]-\omega_{M},\omega_{M}[
  \end{array}
  \right.
\label{eq:4thOMfinal}
\end{small}\end{equation}
and hence MPI sets can be computed according to the methodology presented in the next Section~\ref{sec:Inner}. However, the impact of the model approximation on the MPI set is unknown. Section~\ref{sec:Robust} 
explains how to handle modelling errors.

\section{Inner Approximation of the MPI Set for Polynomial Systems}
\label{sec:Inner}

Let $f$ be a polynomial vector field on $\mathbb{R}^n$. For $x_0 \in \mathbb{R}^n$ and $t\geq 0$, let $x(t|x_0)$ denote the solution of the ordinary differential equation $\dot{x}(t) =f(x(t))$ with initial condition $x(0|x_0)=x_0$. Let $X$ be a bounded and open semi-algebraic 
set of $\mathbb{R}^n$ and $\Bar{X}$ its closure. The Maximal Positively Invariant (MPI) set included in $X$ is defined as:
$$X_0 := \{x_0\in X : \forall t\geq 0,\: x(t|x_0) \in X\}.$$
In words, it is the set of all initial states generating trajectories staying in $X$ \textit{ad infinitum}.

In \cite{korda2013inner}, the authors propose to obtain ROA inner approximations by computing outer approximations of the complementary set (which is a ROA too) with the method presented in \cite{henrion2014convex}. Following the same idea, we chose to focus on approaching the MPI complementary set by the outside in order to get MPI set inner approximations. The MPI complementary set is:
$$ X \backslash X_0 = \{x_0 \in X : \exists t \geq 0, \:x(t|x_0) \in X_\partial\}$$ 
where $X_\partial$ denotes the boundary of $X$. Thus $X \backslash X_0$ is the infinite-time ROA of $X_\partial$. 


This specificity, together with the presence of an algebraic constraint, makes the application of the ROA calculation method as published in the literature not straightforward.  On the one hand, we handle the algebraic constraint by changing the reference measure from an $n$-dimensional volume (Lebesgue) to a uniform measure over a cylinder (Hausdorff). 

On the other hand, ROA approaches usually consider a finite time for reaching the target set. To tackle this issue, we propose here to extend to continuous-time systems the work presented in~\cite{magron2017semidefinite}, where occupation measures are used to formulate the infinite-time reachable set computation problem for discrete-time polynomial systems. 
For a given $a>0$, we define the following linear programming problem:
\begin{equation}
\begin{small}
\begin{array}{rcl}
p^a \:= & \sup & \mu_0(\Bar{X})\\
& \text{s.t.} & \text{div}(f\mu)+\mu_0=\mu_T  \\
&&  \mu_0+\hat{\mu}_0 = \lambda \\
&& \mu(\Bar{X}) \leq a \\
&& 
\end{array}\tag{$P^a$}
\label{eq:Pa}
\end{small}
\end{equation}
where the supremum is with respect to measures $\mu_0 \in \mathcal{M}^+(\Bar{X})$, $\hat{\mu}_0\in \mathcal{M}^+(\Bar{X})$,
$\mu_T \in \mathcal{M}^+(X_\partial)$ and $\mu \in \mathcal{M}^+(\Bar{X})$ with $\mathcal{M}^+(A)$ denoting the cone of non-negative Borel measures supported on the set $A$. 

The first constraint $\text{div}(f\mu)+\mu_0=\mu_T$, a variant of the Liouville equation, encodes the dynamics of the system and ensures that $\mu_0$ (initial measure), $\mu$ (occupation measure) and $\mu_T$ (terminal measure) describe trajectories hitting $X_\partial$. This equation should be understood in the weak sense, i.e. $\forall v \in C^1(\Bar{X}), \int_{\Bar{X}} \text{grad}\: v \cdot f \: d\mu = \int_{X_\partial} v \: d\mu_T - \int_{\Bar{X}} v \: d\mu_0$.

The second constraint ensures that $\mu_0$ is dominated by reference measure: we use a slack measure $\hat{\mu}_0$ and require that $\mu_0 + \hat{\mu}_0 = \lambda$. The third constraint ensures the compacity of the feasible set in the weak-star topology. 

Program (\ref{eq:Pa}) then aims at maximizing the mass of the initial measure $\mu_0$ being dominated by the  reference measure $\lambda$ and supported only on $\Bar{X}\backslash X_0$.  

The dual problem of (\ref{eq:Pa}) is the following linear programming problem: 
\begin{equation}
\begin{small}
\begin{array}{rcl}
d^a \:= & \inf & \displaystyle \int_{\Bar{X}} w(x) \ d\lambda(x)\:+\: u\: a\\
& \text{s.t.} & \text{grad}\: v \cdot f (x) \leq u, \:\forall x \in \Bar{X} \\
&&  w(x) \geq v(x)+1, \forall x \in \Bar{X} \\
&& w(x) \geq 0, \forall x \in \Bar{X} \\
&& v(x) \geq 0, \forall x \in X_\partial
\end{array}\tag{$D^a$}
\label{eq:da}
\end{small}
\end{equation}
where the infimum is with respect to $u \geq 0$, $v \in C^1(\Bar{X})$ and $w \in C^0(\Bar{X})$.

\begin{property}
If 
$(0,v,w)$ is feasible in 
(\ref{eq:da}), then $\{x \in X : v(x)<0\} \subset X_0$ is positively invariant.
\label{lem_pos_inv}
\end{property}

The proof of this statement follows as in \cite[Lemma 2]{henrion2014convex} by evaluating the inequalities in (\ref{eq:da}) on a trajectory. 
Hence, any feasible solution $(0,v,w)$ provides a positively invariant set, thus an inner approximation of the MPI set.

In the same manner than \cite{henrion2014convex,korda2013inner,MCIOUTER}, we use the Lasserre SDP moment relaxation hierarchy of  (\ref{eq:Pa}), denoted $(P^a_k)$, to approach its optimum, where $k \in \mathbb{N}$ is the relaxation order. For brevity and practical reasons (see Fig.~\ref{fig:algo}) this paper presents only the dual hierarchy of SDP SOS tightenings of (\ref{eq:da}): 
\begin{equation}
\begin{small}
\begin{array}{rcl}
d_k^a \:= & \inf & w'l + u a \\
& \text{s.t.} & u - \text{grad}\: v \cdot f = q_{0} + \sum_i q_i\: g_i \\
& & w - v - 1 = p_{0} + \sum_i p_i \: g_i \\
& & w = s_0 + \sum_i s_i \: g_i \\
& & v = t_0 + \sum_i t^+_i g_i - \sum_i t^-_i g_i
\end{array}\tag{$D^a_k$}
\label{eq:dka}
\end{small} 
\end{equation}
where the infimum is with respect to $u \geq 0$, $v \in \mathbb{R}_{2k}[x]$, $w\in \mathbb{R}_{2k}[x]$ and $q_i,p_i,s_i,t_i^+,t_i^- \in \Sigma[x]$, $i=0,1,\ldots,n_X$ with $\mathbb{R}_{2k}[x]$ denoting the vector space of real multivariate polynomials of total degree less than or equal to 2$k$ and $\Sigma[x]$ denoting the cone of SOS polynomials.

For a sufficiently large value of $a>0$ - typically greater than the average escape time on $X\backslash X_0$ - the  optimal solution $(u_k,v_k,w_k)$ of SDP problem (\ref{eq:dka}) is such that $u_k = 0$. 
Hence, solving the SDP program  (\ref{eq:dka}) provides $X_{0,k} := \{x \in X : v_k(x)<0\}$ that is positively invariant from Property \ref{lem_pos_inv}. 
Thus, $X_{0,k}$ is guaranteed to be an \textbf{inner approximation of $X_0$}. 

The algorithmic complexity of the method is that of solving an SDP program whose size is in $O(\binom{n+2k}{2k})$, hence polynomial in the relaxation order $k$ with exponent the number of variables $n$. 

We can now compute inner approximations of the MPI sets of polynomial systems constrained to a semi-algebraic set. However, as discussed before, SM models for transient stability analysis are not polynomial. Although we reformulated them, truncation error may destroy the conservativeness guarantee provided by the proposed method. 

Nevertheless, modelling errors can be seen as an uncertain parameter $\epsilon \in \mathcal{B}$. Hence, there is a compact set $\mathcal{B}\subset \mathbb{R}^p$ such that the non-polynomial vector field $g$ of $\mathbb{R}^n$ satisfies: $$\forall x \in X, \exists \epsilon(x) \in \mathcal{B}, g(x) = f(x,\epsilon(x))$$
where $f$ is a polynomial function from $\mathbb{R}^{n+p}$ to $\mathbb{R}^n$. For instance, the p\textsuperscript{th} order Taylor expansion of $\frac{1}{1+\omega}$ gives: $$\frac{1}{1+\omega} = 1 - \omega +... + (-\omega)^p + \epsilon_p(\omega)$$ with $\epsilon_p(\omega) := \frac{(-\omega)^{p+1}}{1+\omega}$. Thus, $|\epsilon_p(\omega)|\leq \frac{\omega_M^{p+1}}{1-\omega_M}$. In the next section, we propose a robust formulation of the MPI set calculation that ensures the conservative nature of the solution in spite of the modelling errors described by the set $\mathcal{B}$.

\section{Robust MPI Sets}
\label{sec:Robust}

We assume now that the dynamic system depends also on an uncertain time-varying parameter $\epsilon$ evolving in a compact set $\mathcal{B} \subset \mathbb{R}^p$. We are now studying the following ordinary differential equation: $$\dot{x}(t) = f(x(t),\epsilon(t))$$
whose solution is now denoted $x(t|x_0,\epsilon)$ to emphasize the dependence on both the initial condition $x_0$ and the uncertain parameter $\epsilon$. Accordingly, we define the Robust Maximal Positively Invariant (RMPI) set $X_\mathcal{B}$ included in $X$:
$$X_\mathcal{B}:= \{x_0 \in X : \forall \epsilon \in \mathcal{L}^\infty(\mathbb{R}_+,\mathcal{B}),\: \forall t\geq 0, x(t|x_0,\epsilon) \in X \}$$
where $\mathcal{L}^\infty(\mathbb{R}_+,\mathcal{B})$ denotes the vector space of essentially bounded functions from $\mathbb{R}_+$ to $\mathcal{B}$. 
If the system is initialized in $X_\mathcal{B}$, it cannot be brought out of set $X$ by any (time-varying) control whose values belong to $\mathcal{B}$.
Moreover, $X_\mathcal{B}$ is the biggest set included in $X$ being positively invariant for every dynamical system $\dot{x} = f(x,\epsilon)$ with a fixed $\epsilon \in \mathcal{B}$.

In order to compute the RMPI set, we propose the following linear programming problem:
\begin{equation}
\begin{small}
\begin{array}{rcl}
p^a_{\mathcal{B}} \: = & \text{sup} & \mu_0(\Bar{X})\\
&  \text{s.t.} & \text{div}(f\mu)+\mu_T=\mu_0 \\
& & \mu_0+\hat{\mu}_0 = \lambda \\
& & \mu(\Bar{X} \times \mathcal{B}) \leq a
\end{array}\tag{$P^a_\mathcal{B}$}
\label{eq:pab}
\end{small}
\end{equation}
where the supremum is with respect to $\mu_0 \in \mathcal{M}^+(\Bar{X})$, $\hat{\mu}_0\in \mathcal{M}^+(\Bar{X})$,
$\mu_T \in \mathcal{M}^+(X_\partial)$ and $\mu \in \mathcal{M}^+(\Bar{X} \times \mathcal{B})$.

Its dual linear program reads:
\begin{equation}
\begin{small}
\begin{array}{rcl}
d^a_{\mathcal{B}} \:= & \inf & \displaystyle \int_{\Bar{X}} w(x) \ d\lambda(x) \: + \: u \: a\\
&  \text{s.t.} & \text{grad}\: v(x) \cdot f (x,\epsilon) \leq u, \: \forall (x,\epsilon) \in \Bar{X} \times \mathcal{B} \\
& & w(x) \geq v(x)+1, \: \forall x \in \Bar{X} \\
& & w(x) \geq 0, \: \forall x \in \Bar{X} \\
& & v(x) \geq 0, \: \forall x \in X_\partial\\
\end{array}\tag{$D^a _\mathcal{B}$}
\label{eq:dab}
\end{small}
\end{equation}
where the infimum is with respect to $u\geq 0$, $v \in C^1(\Bar{X})$ and $w \in C^0(\Bar{X})$.

\begin{property}
If 
$(0,v,w)$ is feasible  in (\ref{eq:dab}), then the set $\{x \in X : v(x)<0\}$ is positively invariant for any 
given $\epsilon \in \mathcal{L}^\infty(\mathbb{R}_+,\mathcal{B})$.
\label{lem_pos_inv_robuste}
\end{property}

Such a feasible solution is obtained following the same approach as in Section~\ref{sec:Inner}, computing the Lasserre moment hierarchy of (\ref{eq:pab}). Indeed, the dual hierarchy is made of SOS tightenings of (\ref{eq:dab}), which can be solved using SDP.

\section{Numerical Results}
\label{sec:Results}

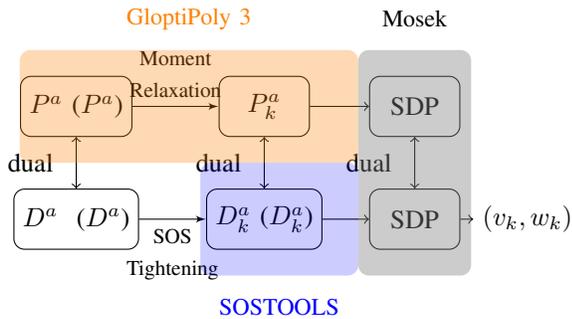
\begin{figure}[H]
\centering
\begin{tikzpicture} [minimum width =  1.2cm,minimum height =  0.8cm]
\node[block]  (Pa)   {$P^a$~\eqref{eq:Pa}};
\node[block,below of=Pa,node distance=1.5cm]  (Da)   {$D^a$ ~\eqref{eq:da}};
\node[block,right of=Pa,node distance=2.5cm]  (Pak)   {$P^a_k$};
\node[block,right of=Da,node distance=2.5cm]  (Dak)   {$D^a_k$~\eqref{eq:dka}};
\node[block,right of=Pak,node distance=2cm]  (SDPp)   {SDP};
\node[block,right of=Dak,node distance=2cm]  (SDPd)   {SDP};
\node[right of=SDPd,node distance=1.5cm]  (poly)   {$(v_k,w_k)$};
\path [line] (Pa) -- node [text width=1.5cm,midway,above, text centered] {\footnotesize Moment Relaxation} (Pak);
\path [line] (Da) -- node [text width=1.5cm,midway,below, text centered] {\footnotesize SOS Tightening} (Dak);
\draw[<->] (Pa) -- node[left]{dual} (Da);
\draw[->] (Pak) -- (SDPp);
\draw[<->] (Pak) -- node[left]{dual} (Dak);
\draw[<->] (SDPp) -- node[left]{dual} (SDPd);
\draw[->] (Dak) -- (SDPd);
\draw[->] (SDPd) -- (poly);
\node[fill=gray,rounded corners,below of=SDPp,node distance=0.75cm,minimum height = 3cm, minimum width =  1.5cm, opacity=.4,label={\small Mosek}] {};
\node[fill=orange,rounded corners,right of=Pa,node distance=1.5cm,minimum height = 1.5cm, minimum width =  4.5cm, opacity=.3,label={\small \textcolor{orange}{GloptiPoly 3}}] {};
\node[fill=blue,rounded corners,right of=Dak,node distance=0.2cm,minimum height = 1.5cm, minimum width =  2.1cm, opacity=.2,label=below:{\small \textcolor{blue}{SOSTOOLS}}] {};
\end{tikzpicture}
\vspace{-2mm}
\caption{Implementation of the proposed method}
\label{fig:algo}
\end{figure}

The method described in this work has been implemented in MATLAB. The SDP problems are solved using MOSEK that takes as input a raw 
SDP program. As illustrated in Fig.~\ref{fig:algo} 
we consider two equivalent alternatives to produce this file:  

\begin{enumerate}
    \item Using the interface GloptiPoly 3 \cite{henrion2009gloptipoly}, that takes a linear program on measures as an input and produces the Lasserre SDP moment relaxation of a specified degree.
    \item Using the interface SOSTOOLS \cite{prajna2002introducing}, that takes a SOS programming problem as an input and produces the corresponding SDP problem.
\end{enumerate}

\vspace{-5mm}

\vspace{-2mm}

It is important to highlight that from the implementation point of view, the SM models presented in Section~\ref{sec:Modelling} were renormalized in order to get well scaled SDP problems. In addition, a \textit{reasonable} set $X$ is defined such that the volume of the MPI set $X_0$ covers a non-negligible part of this box. For the test systems considered here, this was achieved by setting all variables between -1 and 1. For parameter $a$, we used 100.

\subsection{Link between MPI sets and transient stability}

\begin{figure}[h!]
\centering
\begin{subfigure}{0.25\textwidth}
\centering
    \includegraphics[width=\textwidth,trim={0 0 2cm 15cm},clip]{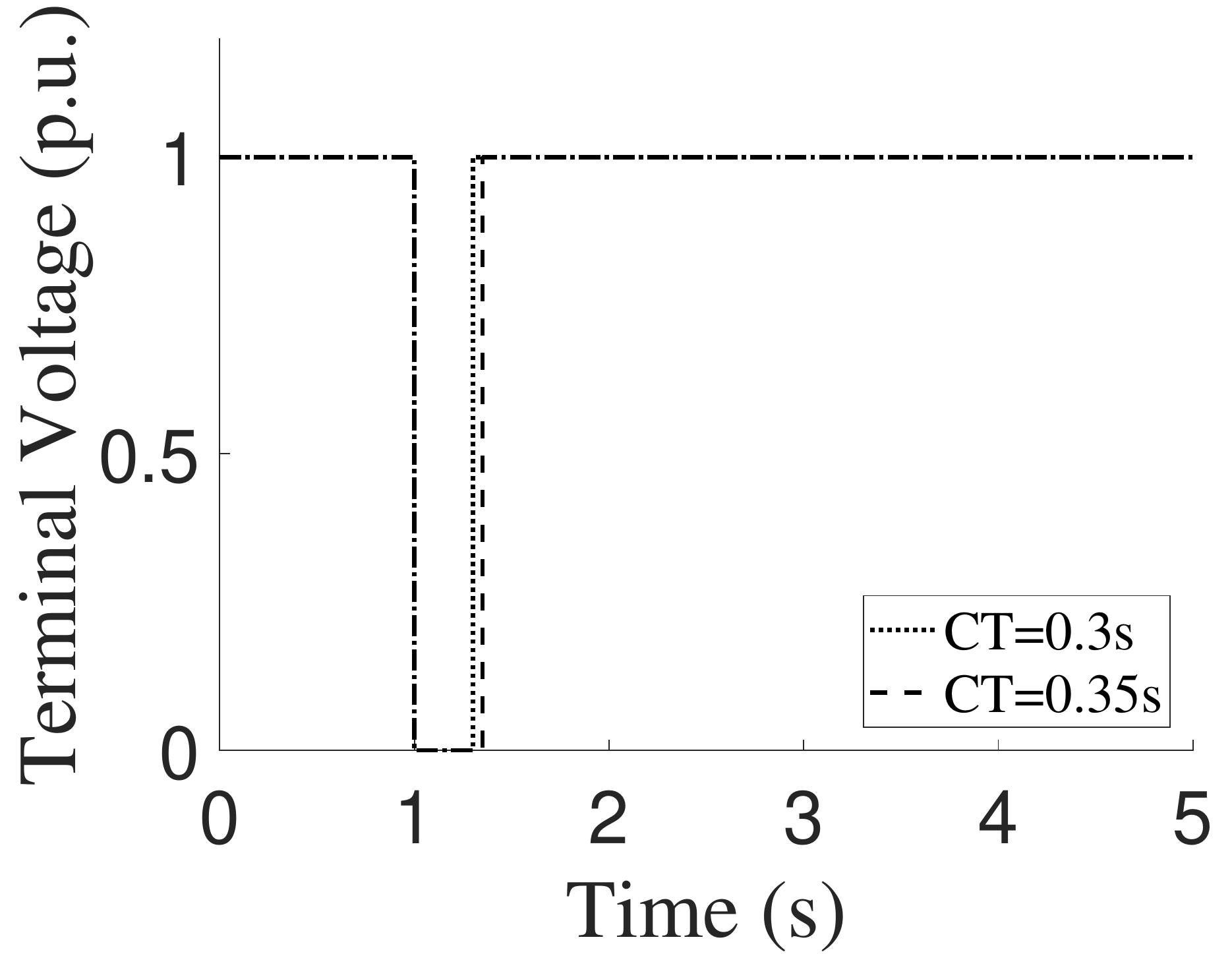}
    \caption{Terminal voltage}
    \label{fig:Vs}
    \end{subfigure}%
    \begin{subfigure}{0.25\textwidth}
    \centering
    \includegraphics[width=\textwidth,trim={0 0 2cm 15cm},clip]{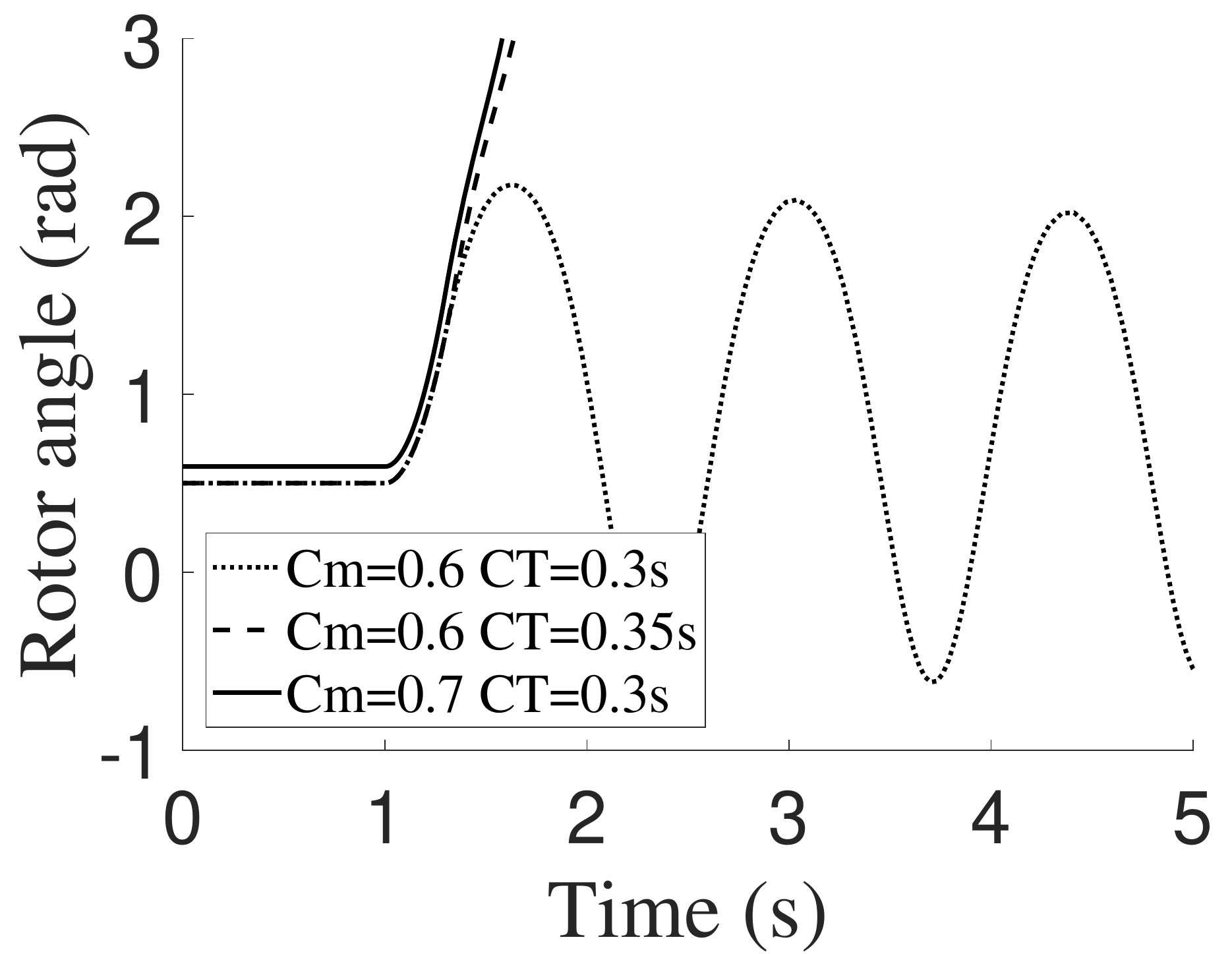}
    \caption{Rotor angle trajectory}
    \label{fig:RotAng}
    \end{subfigure}
\caption{Stable and unstable cases}
\label{fig:2ndOMVsDelta}
\end{figure}

Let us consider: i) the test system described in Appendix~\ref{sec:appB}, ii) two scenarios with $C_m=$0.6 p.u. and $C_m=$0.7 p.u., and iii) for illustrative purposes, two faults at the SM terminal with different clearing times (CT): 300 and 350 ms, see Fig.~\ref{fig:Vs}. 

\subsubsection{SM 2\textsuperscript{nd} OM}
Critical clearing times for both scenarios are determined through simulation. For this first model $CCT_1$ for scenario 1 ($C_m=$0.6 p.u) is 310 ms and $CCT_2$ for scenario 2 ($C_m=$0.7 p.u) is 250 ms. Figure~\ref{fig:RotAng} shows that when the fault is cleared at 300 s only scenario 1 remains stable.

\begin{figure}[h!]
\centering
\includegraphics[width=0.5\textwidth,trim={0 0 1cm 20cm},clip]{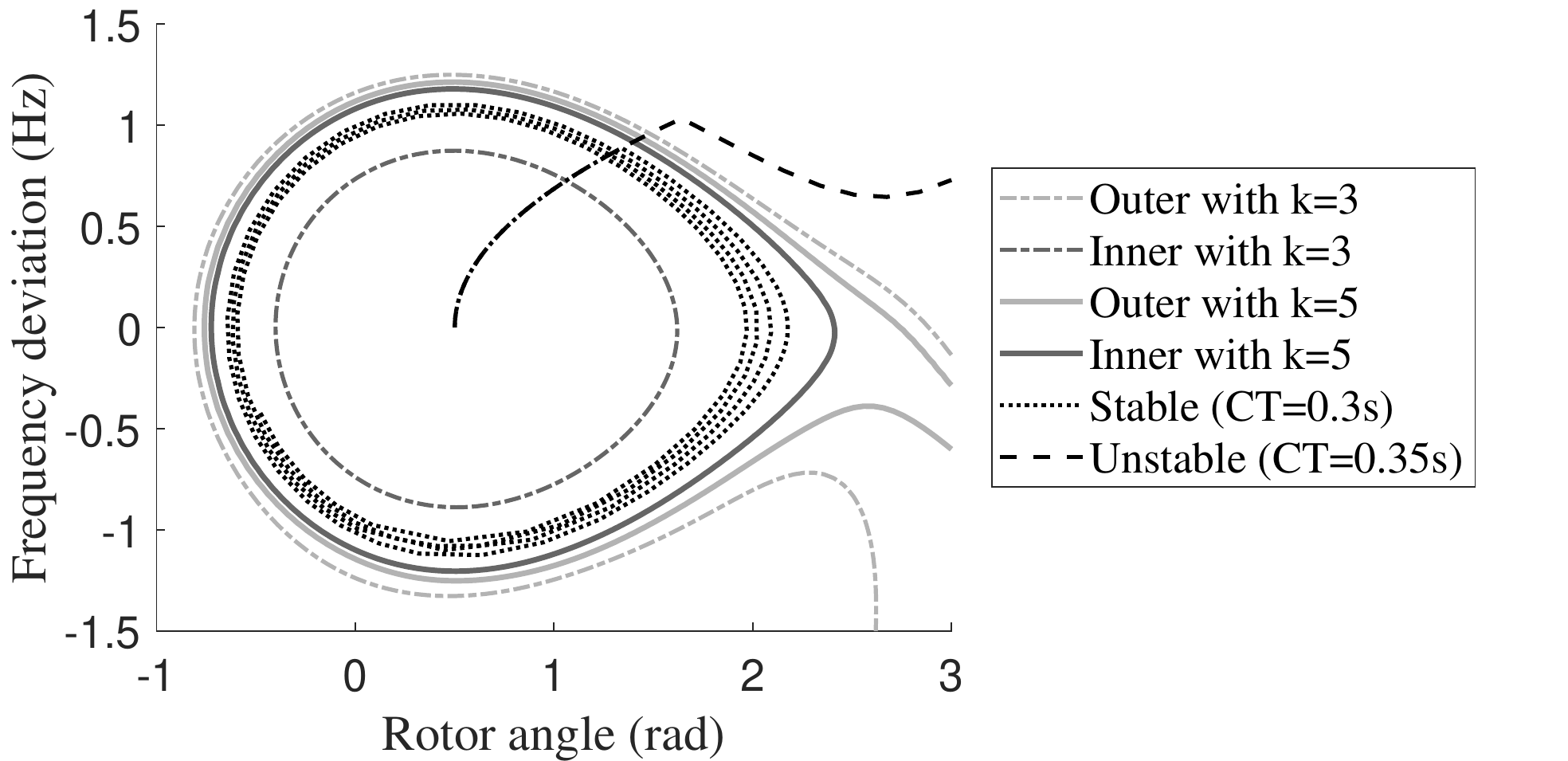}
\caption{MPI set approximations for different $k$.}
\label{fig:2ndInOutK}
\end{figure}

Figure~\ref{fig:2ndInOutK} shows the MPI set approximations computed with the proposed approach -solving program~\eqref{eq:dka} and setting $v=$0- for two different degrees of relaxation ($k=3$ and $k=5$). 
The trajectories presented in Fig.~\ref{fig:RotAng} for scenario 1 and different fault clearing times are also included. The accuracy gain provided by increasing the relaxation degree is observed. 

For $k=$5 the computed inner and outer approximations are quite close, which enables us to conclude about the stability of a certain post-fault situation by simulating only the faulted system. Moreover, they provide an insight on the "stability margin" by looking into the distance between the system state at the fault elimination and the boundaries of the MPI set. 

\begin{figure}[h!]
\centering
\includegraphics[width=0.45\textwidth,trim={0 0 1cm 20cm},clip]{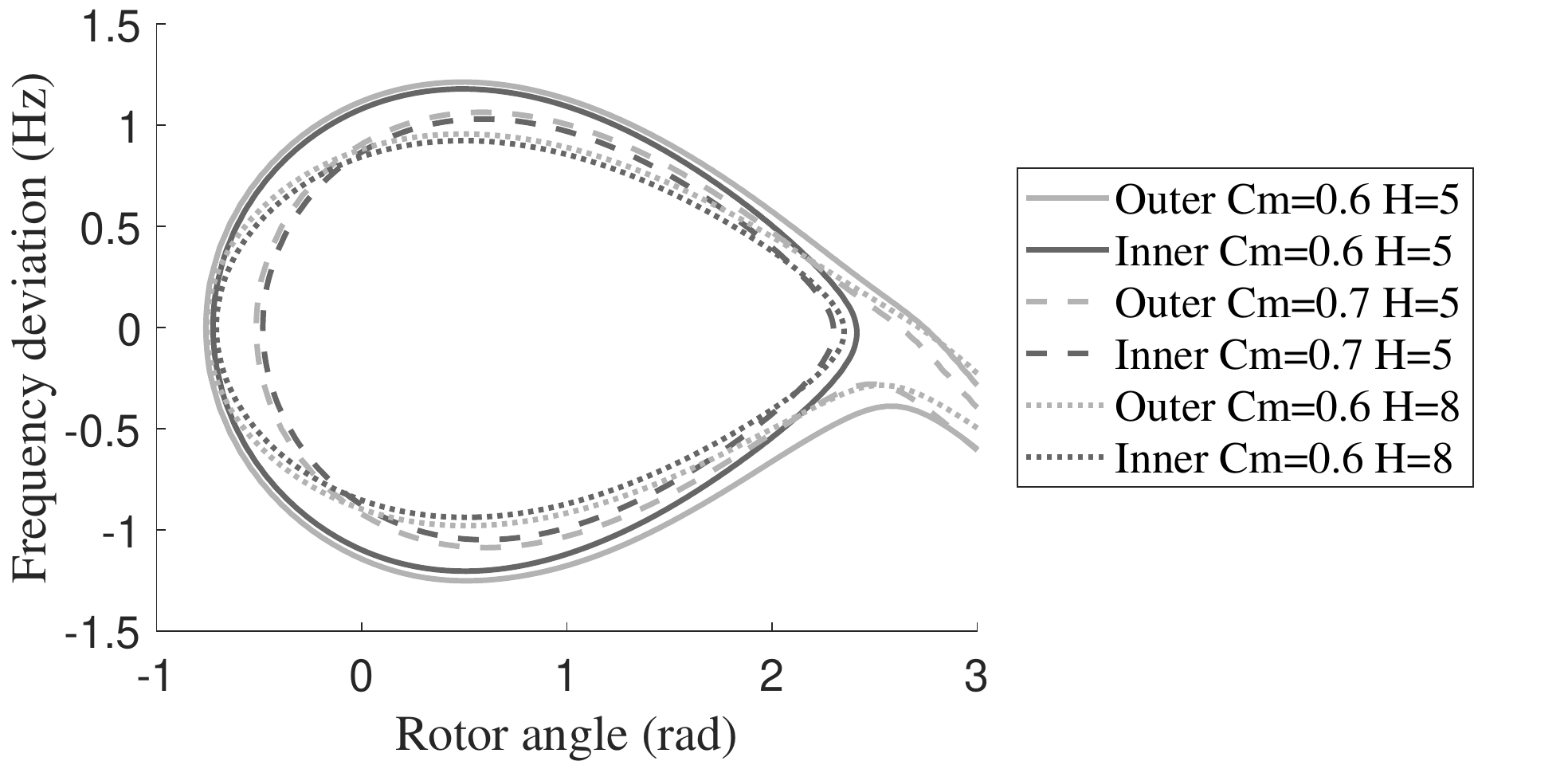}
\caption{MPI set approximations for different parameters.}
\label{fig:2ndCmH}
\vspace{-5mm}
\end{figure}

Figure~\ref{fig:2ndCmH} shows the computed MPI set for $k=5$ and different system parameters. In all cases CPU times are around 4 seconds\footnote{Intel(R) Core(TM) i7-4900MQ 2.8GHz.}. As the SM is operated closer to its maximal capacity (higher $C_m$) the stability region becomes smaller and moves to the right side. Indeed, the rotor angle at the equilibrium point increases with $C_m$.   

Consistently with intuitions from the equal area criterion, the critical angle after which the fault elimination becomes ineffective to prevent loss of synchronism is independent of $H$. However, the MPI sets become "flatter" because the lower the inertia of the unit, the higher the speed that will be possible to arrest for the same available decelerating power. 

\begin{figure}[h!]
\centering
\begin{subfigure}{0.25\textwidth}
\centering
    \includegraphics[width=\textwidth,trim={0 0 2cm 15cm},clip]{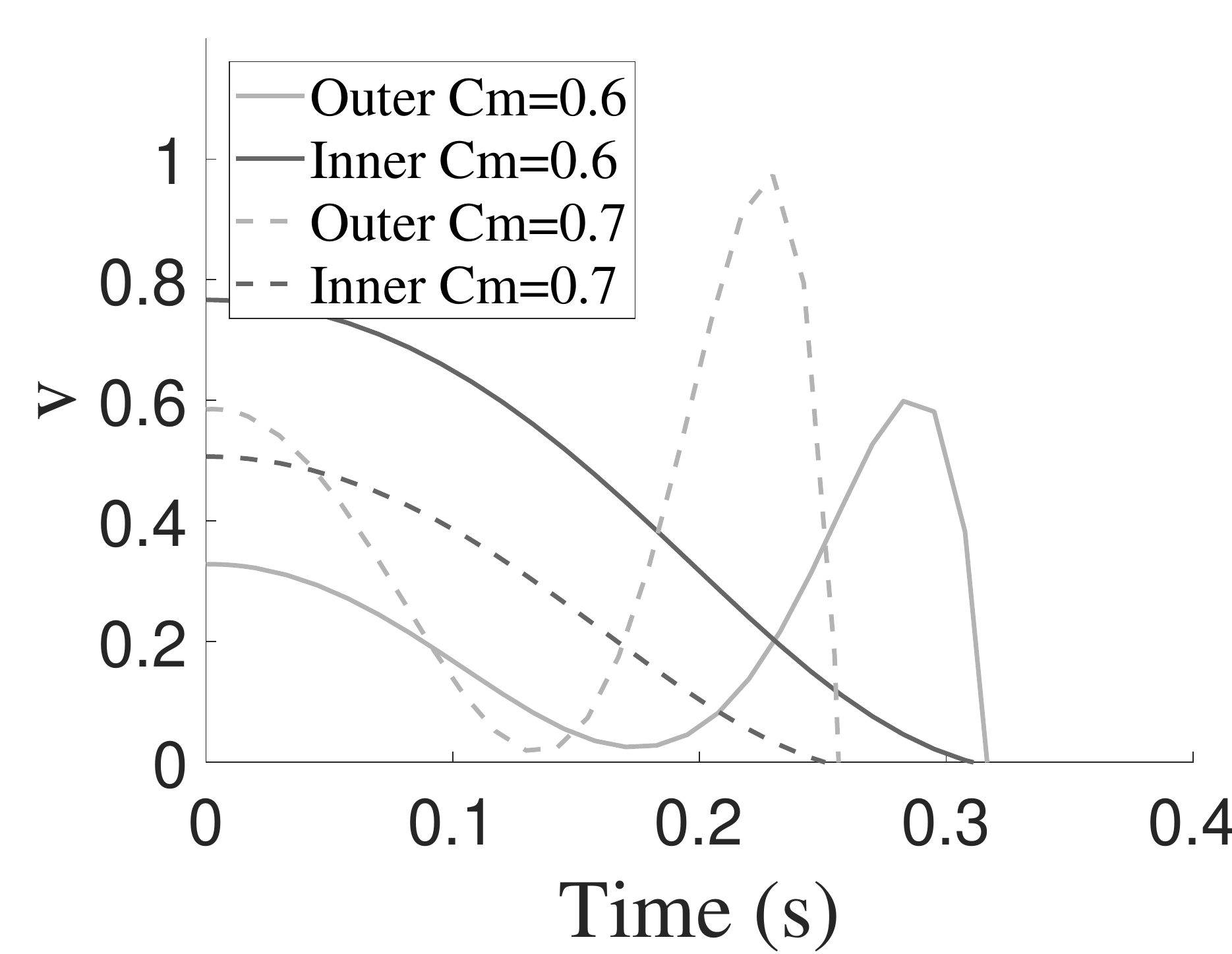}
    \caption{CCT for $H=$5 MWs/MVA}
    \label{fig:CCTH1}
    \end{subfigure}%
    \begin{subfigure}{0.25\textwidth}
    \centering
    \includegraphics[width=\textwidth,trim={0 0 2cm 15cm},clip]{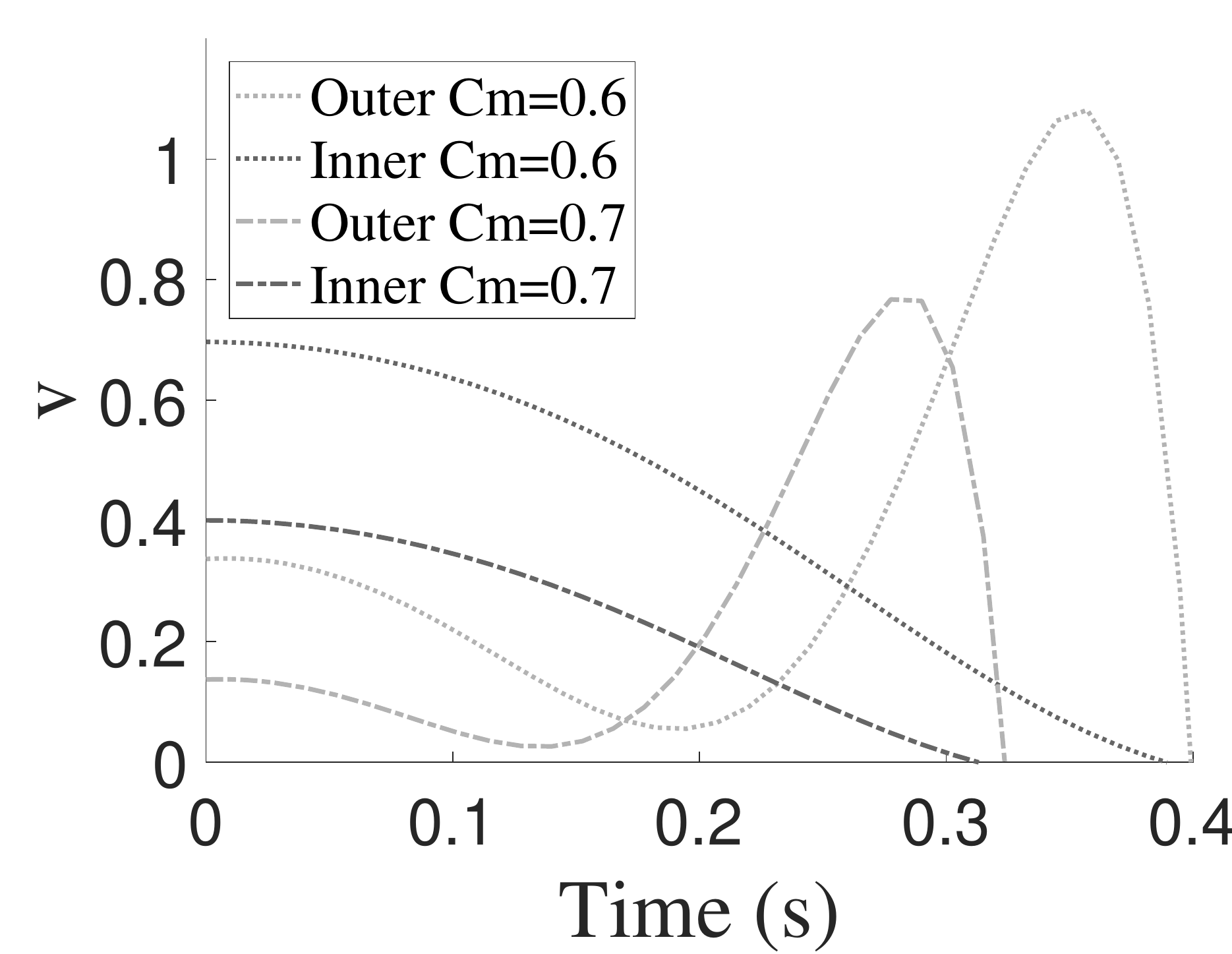}
    \caption{CCT for $H=$8 MWs/MVA}
    \label{fig:CCTH2}
    \end{subfigure}
\caption{$v$ evaluated on the faulted trajectory}
\label{fig:2ndCCT}
\end{figure}

Of course, as shown in Fig.~\ref{fig:2ndCCT}, the CCT for a given fault increases with $H$. These figures show the polynomial $v(\delta,\omega)$, describing the MPI set and obtained by solving~\eqref{eq:dka}, evaluated along the faulted trajectory. 

\vspace{-3mm}

\vspace{-2mm}

For readability purposes the sign of $v$ for the inner approximation has been changed and  the values have been normalized. The zero crossing with the abscissa axis corresponds to the moment when $v(\delta(t),\omega(t)) = 0$, which means that the state variables are no longer inside the computed MPI set. This is the CCT. For $k=5$ the CCT is estimated with a 10 ms precision. However, $v$ remains a high degree polynomial.

\subsubsection{SM 3\textsuperscript{rd} OM}
\begin{figure}[h!]
\centering
\captionsetup[subfigure]{labelformat=empty}
\begin{subfigure}{0.25\textwidth}
\centering\includegraphics[width=\textwidth,trim={0 2.75cm 2cm 11.75cm},clip]{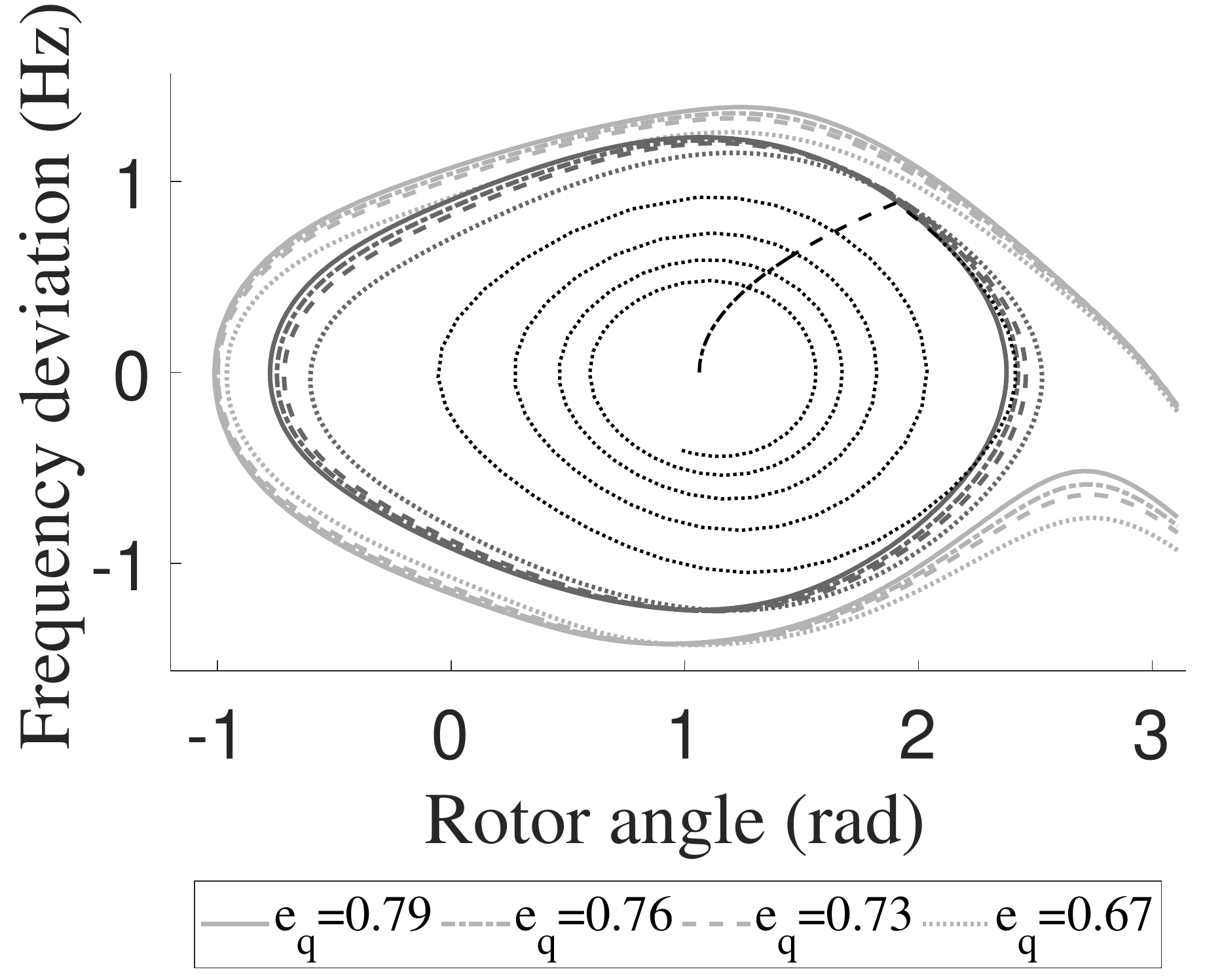}
    \end{subfigure}%
    \begin{subfigure}{0.25\textwidth}
    \centering
    \includegraphics[width=\textwidth,trim={0 2.75cm 2cm 11.75cm},clip]{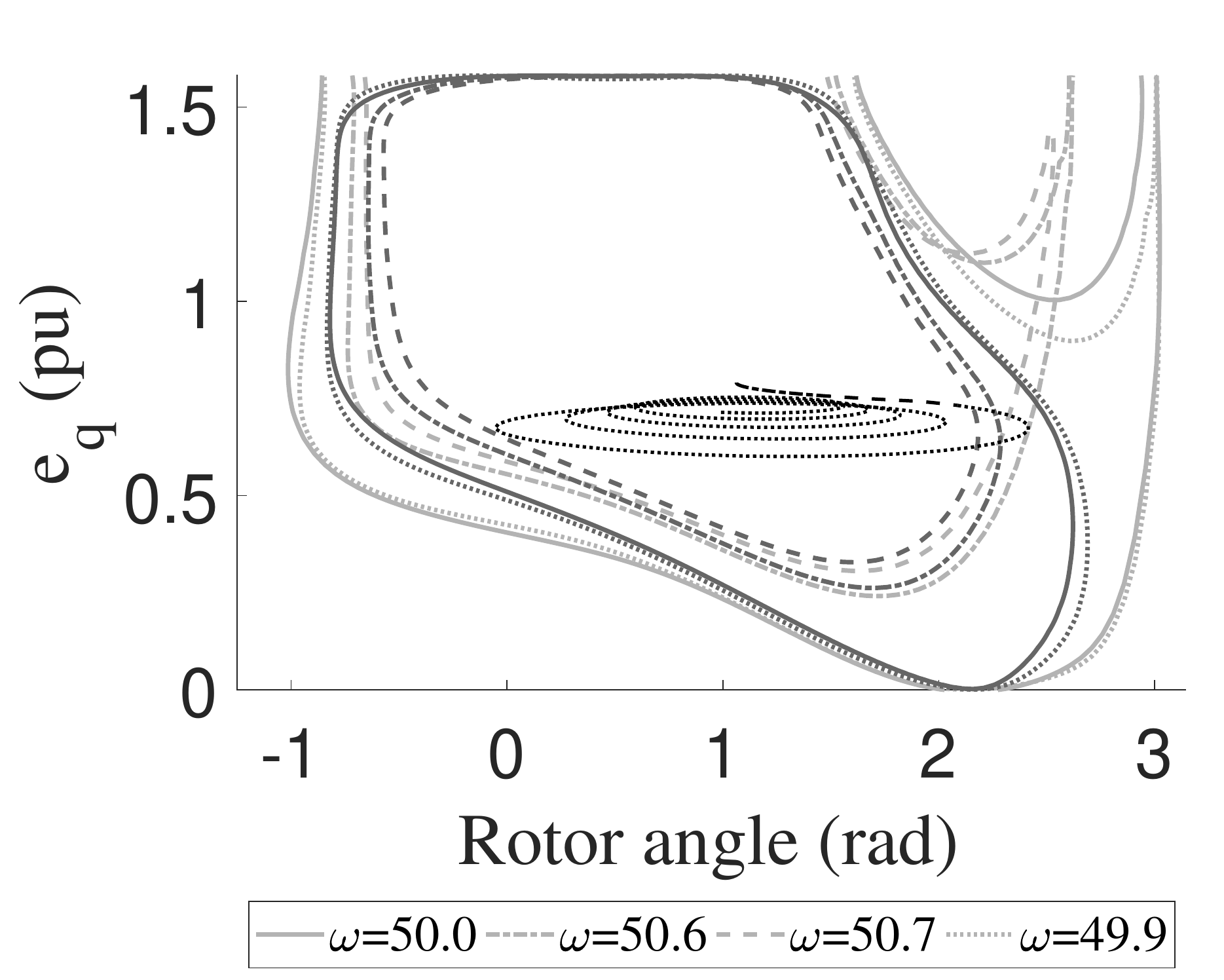}
    \end{subfigure}
\caption{Sections of the MPI set for the 3\textsuperscript{rd} OM}
\label{fig:MPI3rdOM}
\vspace{-4mm}
\end{figure}

In this case, the MPI set consists in a three dimensional volume ($\delta$,$\omega$,$e'_q$). Figure~\ref{fig:MPI3rdOM} shows sections of the outer (light grey) and inner (dark grey) approximations of the MPI set for different values $e'_q$ (left) and $\omega$ (right). These values correspond to specific points of the stable fault trajectory ($t=$1s, $t=$1.2s, $t=$1.35s and $t=$1.5s respectively). As expected, the stability region is larger for high values of internal electromotive forces (the set $X$ is limited to $2e'_{q0}=$1.58pu). Again, evaluating $v$ along the trajectory during the fault enables us to bound the CCT between 305-330ms.

\subsubsection{SM 4\textsuperscript{th} OM} 
\begin{figure}[h!]
\centering
\begin{subfigure}{0.25\textwidth}
\centering\includegraphics[width=\textwidth,trim={0 2.75cm 2cm 15cm},clip]{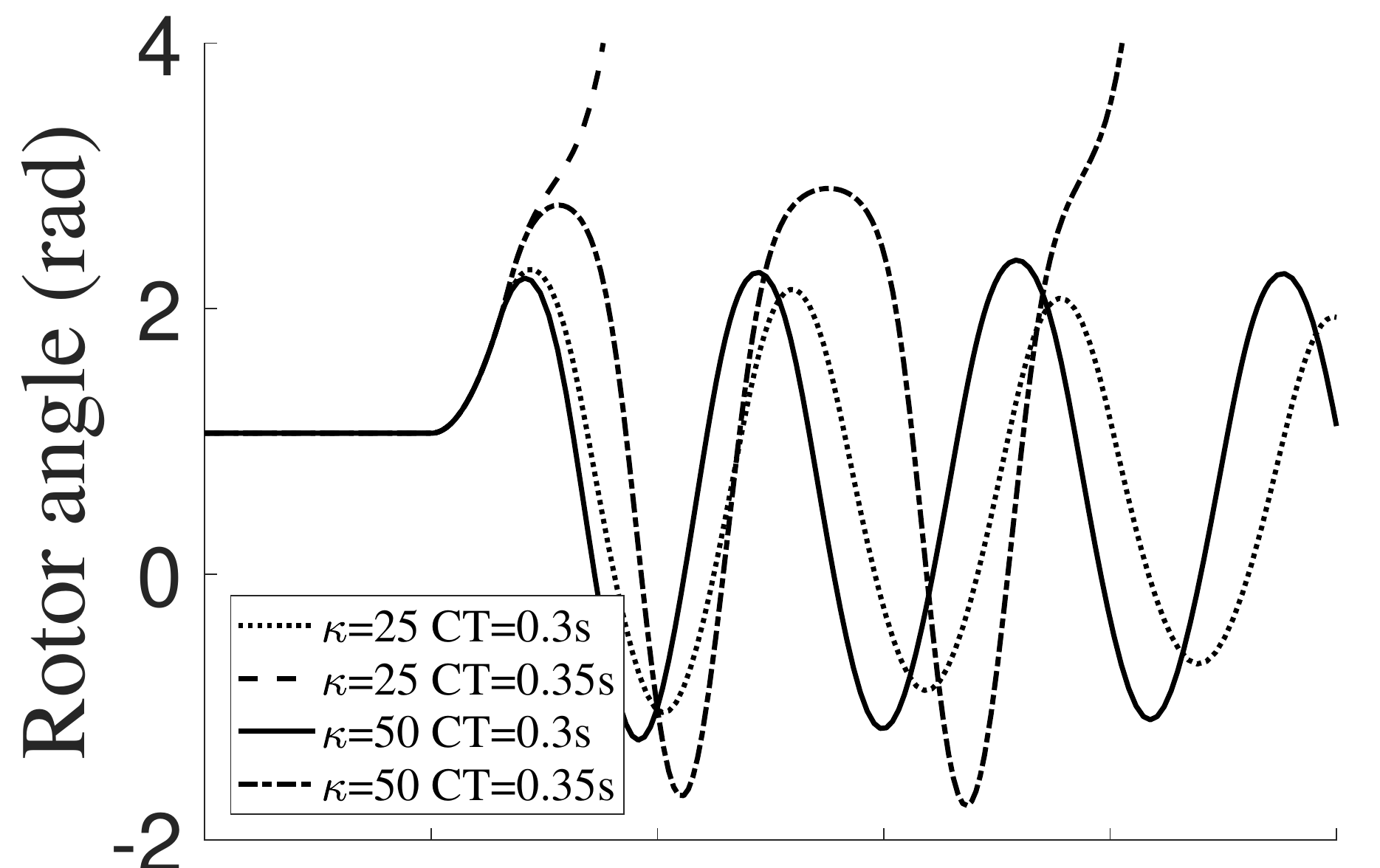}
\caption{Trajectory during faults 4\textsuperscript{th} OM}
\label{fig:4thOMrotang}
    \end{subfigure}%
    \begin{subfigure}{0.25\textwidth}
    \centering
    \includegraphics[width=\textwidth,trim={0 2.75cm 2cm 15cm},clip]{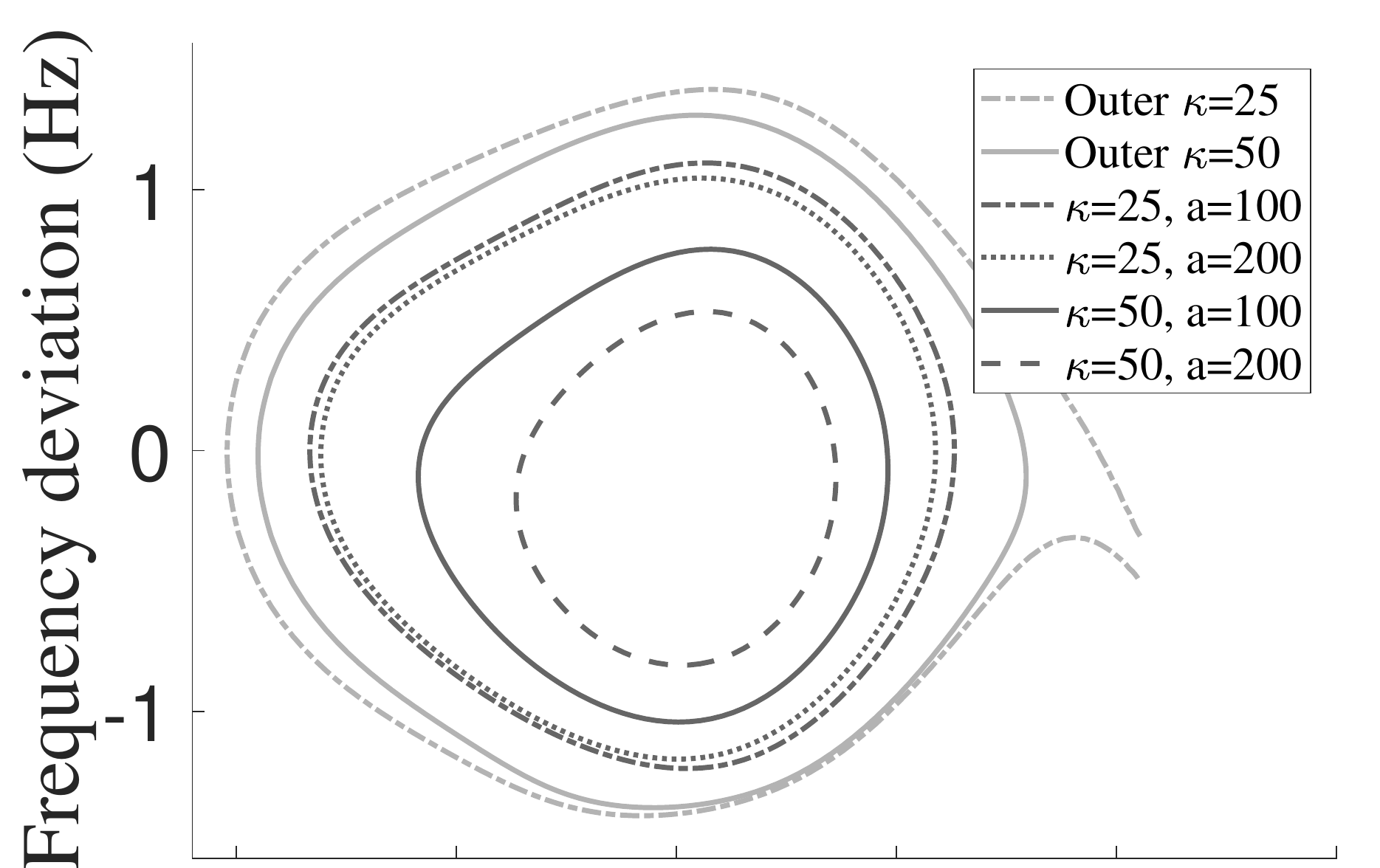}
    \caption{MPI set approximations}
    \label{fig:MPI4thOM}
    \end{subfigure}
\caption{Results for the 4\textsuperscript{th} OM}
\label{fig:Results4thOM}
\vspace{-2mm}
\end{figure}

it is well known that voltage regulators may introduce negative damping in the system~\cite{kundur1994power}. Although power plants have more sophisticated controller, we consider here a proportional one as described in~\eqref{eq:4thOM} for illustrative purposes. Figure~\ref{fig:4thOMrotang} shows that depending on the value of $\kappa$ the lost of synchronism may occurs after a few diverging oscillations. 

Fig~\ref{fig:MPI4thOM} shows sections of the MPI set outer (light grey) and inner (dark grey) approximations at  the equilibrium point. For badly-damping system the choice of parameter $a$ may have a impact in accuracy for a given relaxation degree $k$.

\vspace{-2mm}

\subsection{Model approximation and Robust MPI sets}
\begin{figure}[h!]
\centering\includegraphics[width=0.45\textwidth,trim={0 0 2cm 20cm}]{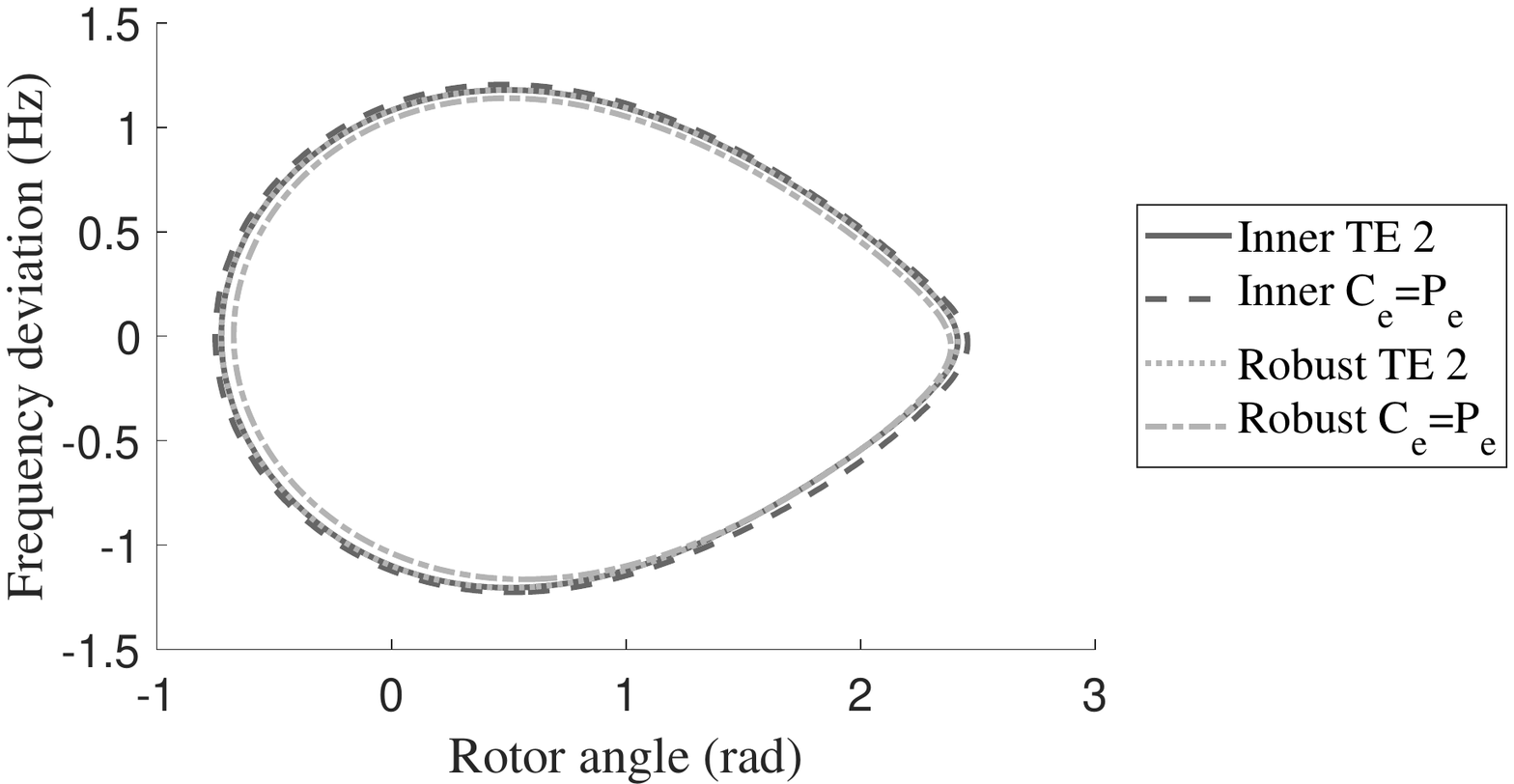}
\caption{Robust MPI set}
\vspace{-2mm}
\label{fig:robust}
\end{figure}

Previous section considered a 2\textsuperscript{nd} order Taylor expansion for the $\frac{1}{\omega_s}$ term as described in Section~\ref{sec:Modelling}. Figure~\ref{fig:robust} shows that this approximation is accurate enough since the RMPI set overlaps. However, in the presence of larger modeling errors, for instance, if we use the electrical power directly into the speed equation ($C_e=P_e$ and $\frac{1}{\omega_s}\approx 1$), we observe  that: 
\begin{enumerate}
    \item The MPI set computed with~\eqref{eq:dka} is wider and may include points that are unstable in the original non-polynomial system.
    \item Since the bounds of the $\mathcal{B}$ are larger, the RMPI set is a bit smaller, but offers conservativeness guarantees. 
\end{enumerate}

\vspace{-2mm}

Indeed, in the second case we write $\frac{1}{1+\omega} = 1 + \epsilon_0(\omega)$ with $\epsilon_0(\omega) = \frac{- \omega}{1 + \omega}$ whereas in the case of the 2\textsuperscript{nd} order Taylor expansion, we write $\frac{1}{1+\omega} = 1 - \omega + \omega^2 + \epsilon_2(\omega)$, with $\epsilon_2(\omega) = \frac{-\omega^3}{1+\omega}$. Naturally, the bounds on $\epsilon_2$ are tighter : $|\epsilon_0(\omega)|\leq \frac{\omega_M}{1-\omega_M}$ while $|\epsilon_2(\omega)| \leq \frac{\omega_M^3}{1-\omega_M}$ (typically $\omega_M = 0.05)$. 


\subsection{Performance for the 3\textsuperscript{rd} and 4\textsuperscript{th} OM}
As discussed before, the algorithmic complexity of the method depends strongly on $n$, the number of states. 
Table~\ref{tab:volume} shows the accuracy on the computation of the MPI set inner and outer approximations with the relaxations degree. The volumes are computed using a Monte-Carlo method. Table~\ref{tab:CPU} presents the associated computing time for both models. It is observed that CPU time raises considerably for the 4\textsuperscript{th} OM. 

\begin{table}[h!]
\caption{Volume of the computed MPI sets}
\begin{center}
\begin{tabular}{|c|p{2mm}|p{2.4cm}|p{2.4cm}|}
\hline
Model & $k$ & inner approximation & outer approximation \\
\hline
3\textsuperscript{rd} & 4 \newline 5 \newline 6 & 9.84 \newline 10.37 \newline 10.85 & 17.02 \newline 14.91 \newline 13.97 \\
\hline 
4\textsuperscript{th} & 4 \newline 5 & 12.00 \newline 17.06 & 30.40 \newline 28.12 \\
\hline
\end{tabular}
\label{tab:volume}
\end{center}
\end{table}

\begin{table}[h!]
\caption{CPU time of the MPI set computation (s)}
\begin{center}
\begin{tabular}{|c|p{2mm}|p{2.4cm}|p{2.4cm}|}
\hline
Model & $k$ & inner approximation & outer approximation \\
\hline
3\textsuperscript{rd} & 4 \newline 5 \newline 6 &  12.64 \newline 29.92 \newline 129.04 & 4.50 \newline 20.13 \newline 100.06 \\
\hline 
4\textsuperscript{th} & 4 \newline 5 & 63.82 \newline 573.65 & 41.04 \newline 339.65  \\
\hline
\end{tabular}
\label{tab:CPU}
\end{center}
\end{table}



\vspace{-4mm}







\section{Conclusions}
\label{sec:Conclusions}
The transient stability problem has been formulated as the inner approximation of the MPI set of the polynomial dynamic system. For this purpose, we have first transformed SM machines models into polynomial ones, and then adapted the published work based on occupation measures and Lasserre hierarchy to the infinite-time ROA calculation for continuous systems constrained to an algebraic set. Simulation results showed that we can compute multidimensional stability regions for more complex SM models and that CCT can be accurately bounded evaluating the obtained polynomial for inner and outer MPI set approximations on the faulted trajectory. 

Moreover, we have proposed a robust formulation that provides conservativeness guarantees in the presence of bounded modelling uncertanties. Again accurate results are obtained when taking into account Taylor approximation errors.

However, algorithmic complexity leads to high CPU times as more details were included in the SM model. Future work will focus on limiting the required relaxation degree in order to reduce computational cost and be able to increase the state space dimension.  


\section*{Acknowledgment}
The authors would like to thanks Matteo Tacchi from LAAS-CNRS Toulouse, and Philippe Juston from RTE for the enlightening discussions.


\appendices


\section{Test System}
\label{sec:appB}

\begin{figure}[h!]
\centerline{\includegraphics[width=0.35\textwidth]{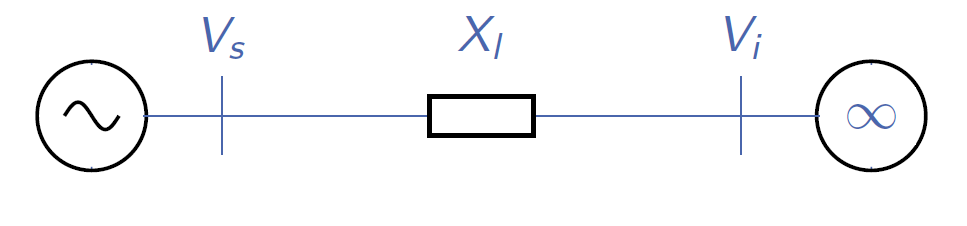}}
\caption{Test system: Single Machine Infinite Bus}
\end{figure}

\begin{table}[h!]
\caption{System Parameters}
\begin{center}
\begin{tabular}{|c|c|c|c|c|}
\hline
$\omega_n$ & 314 $rad.s^{-1}$ & & $\Bar{E}_{vf}$ (3\textsuperscript{rd} OM) & 1.85 pu \\
\hline
$H$ & 5 MWs/MVA & & $T_d$ & 10 s\\
\hline 
$D$ & 1 pu & & $T_E$ & 0.3 s  \\
\hline
$V_s$ & 1 pu & & $\kappa$ & 25 pu  \\
\hline
$V_i$ & 1 pu & & $x_d$ & 2.5 pu\\
\hline
$X_l$ (2\textsuperscript{dn} OM) & 0.8 pu & & $x_q$ & 2.5 pu\\
\hline
$X_l$ (3\textsuperscript{rd} OM) & 0.2 pu & & $x'_d$ & 0.4 pu\\
\hline
\end{tabular}
\label{tab:syspar}
\end{center}
\end{table}



\end{document}